\documentclass[reqno]{amsart}
\usepackage{hyperref}

\usepackage{amsmath,amssymb,amsthm}
\usepackage{enumerate}
\usepackage{color}
\usepackage{epsfig}
\usepackage{graphics}
\usepackage{graphicx}
\usepackage{subfigure}

\begin{document}
\title[]{The Cesaro Operator in $\ell^{2}$ is Essentially Normal}% The title of the article
\author[u. g\"{u}l]{u\u{g}ur g\"{u}l}

\address{u\u{g}ur g\"{u}l,  \newline
Hacettepe University, Department of Mathematics, 06800, Beytepe,
Ankara, TURKEY}
\email{\href{mailto:gulugur@gmail.com}{gulugur@gmail.com}}
%\href{mailto:myaddress@wikibooks.org}{myaddress@wikibooks.org}
%\email{{mailto:aozbekler@gmail.com}{aozbekler@gmail.com}}

\thanks{Submitted XXX}

\subjclass[2000]{47B33} \keywords{Cesaro Operator, Hardy Space of
the disc, Essential Normality.}
\begin{abstract}
In this paper we prove that the Cesaro operator $\mathcal{C}$ in
$\ell^{2}$, the Hilbert space of square summable sequences, is
essentially normal, i.e. the commutator
$[\mathcal{C}^{\ast},\mathcal{C}]:=\mathcal{C}^{\ast}\mathcal{C}-\mathcal{C}\mathcal{C}^{\ast}$
is a compact operator on $\ell^{2}$.
\end{abstract}

\maketitle
%\numberwithin{equation}{section}
\newtheorem{theorem}{Theorem}
\newtheorem{acknowledgement}[theorem]{Acknowledgement}
\newtheorem{algorithm}[theorem]{Algorithm}
\newtheorem{axiom}[theorem]{Axiom}
\newtheorem{case}[theorem]{Case}
\newtheorem{claim}[theorem]{Claim}
\newtheorem{conclusion}[theorem]{Conclusion}
\newtheorem{condition}[theorem]{Condition}
\newtheorem{conjecture}[theorem]{Conjecture}
\newtheorem{corollary}[theorem]{Corollary}
\newtheorem{criterion}[theorem]{Criterion}
\newtheorem{definition}[theorem]{Definition}
\newtheorem{example}[theorem]{Example}
\newtheorem{exercise}[theorem]{Exercise}
\newtheorem{lemma}[theorem]{Lemma}
\newtheorem{notation}[theorem]{Notation}
\newtheorem{problem}[theorem]{Problem}
\newtheorem{proposition}[theorem]{Proposition}
\newtheorem{remark}[theorem]{Remark}
\newtheorem{solution}[theorem]{Solution}
\newtheorem{summary}[theorem]{Summary}
\newtheorem*{thma}{Main Theorem}
\newcommand{\norm}[1]{\left\Vert#1\right\Vert}
\newcommand{\abs}[1]{\left\vert#1\right\vert}
\newcommand{\set}[1]{\left\{#1\right\}}
\newcommand{\Real}{\mathbb R}
\newcommand{\eps}{\varepsilon}
\newcommand{\To}{\longrightarrow}
\newcommand{\BX}{\mathbf{B}(X)}
\newcommand{\A}{\mathcal{A}}
\newtheorem*{schur}{Schur's Boundedness Test}

\section*{introduction}
 The Cesaro operator acting on $\ell^{2}$ is a very well known
 operator and is defined as
 $\mathcal{C}:\ell^{2}\rightarrow\ell^{2}$,
 $$\mathcal{C}\bigg((x_{n})_{n=0}^{+\infty}\bigg):=\bigg(\frac{1}{n+1}\sum_{k=0}^{n}x_{k}\bigg)_{n=0}^{+\infty}.$$
 More precisely $\mathcal{C}$ is defined as
$$\mathcal{C}\bigg(x_{0},x_{1},x_{2},...,x_{n},...\bigg)=\bigg(x_{0},\frac{x_{0}+x_{1}}{2},\frac{x_{0}+x_{1}+x_{2}}{3},...,\frac{1}{n+1}\sum_{k=0}^{n}x_{k},...\bigg).$$
  The Cesaro operator has the following matrix representation on
  $\ell^{2}$:
$$[\mathcal{C}]:=
\begin{pmatrix}
1&0&0&0\dots\\
\frac{1}{2}&\frac{1}{2}&0&0\dots\\
\frac{1}{3}&\frac{1}{3}&\frac{1}{3}&0\dots\\
\vdots & \vdots & \vdots &\vdots\\
\frac{1}{n}&\frac{1}{n}&\frac{1}{n}&\frac{1}{n}\dots &\frac{1}{n}&0&\dots\\
\vdots & \vdots & \vdots &\vdots
\end{pmatrix}
$$
where
$\ell^{2}:=\{(x_{n})_{n=0}^{+\infty}:x_{n}\in\mathbb{C},\sum_{n=0}^{+\infty}\mid
x_{n}\mid^{2}<+\infty\}$ is the Hilbert space of square summable
complex sequences. The adjoint operator $\mathcal{C}^{\ast}$,
then, has the matrix representation:
$$[\mathcal{C}^{\ast}]:=
\begin{pmatrix}
1&\frac{1}{2}&\frac{1}{3}&\frac{1}{4}\dots\\
0&\frac{1}{2}&\frac{1}{3}&\frac{1}{4}\dots\\
0&0&\frac{1}{3}&\frac{1}{4}\dots\\
0&0&0&\frac{1}{4}\dots\\
\vdots & \vdots & \vdots &\vdots\\
0&0&0&0\dots &0&\frac{1}{n}&\dots\\
\vdots & \vdots & \vdots &\vdots
\end{pmatrix}.$$

Let $X$, $Y$ be two Banach spaces over the field of complex
numbers $\mathbb{C}$ and $K:X\rightarrow Y$ be a bounded linear
operator. If the image $K(\mathbb{B}_{X})\subset Y$ of the closed
unit ball $\mathbb{B}_{X}:=\{x\in X:\parallel x\parallel_{X}\leq
1\}$ under $K$ is a norm compact subset of $Y$ then $K$ is called
a \textbf{compact operator}. An bounded linear operator
$T:H\rightarrow H$ acting on a Hilbert space $H$ over the field of
complex numbers $\mathbb{C}$ is called ``\textbf{essentially
normal}" if its self commutator
$[T^{\ast},T]:=T^{\ast}T-TT^{\ast}$ is a compact operator. Here
$T^{\ast}$ denotes the adjoint of $T$ which is the unique operator
$T^{\ast}:H\rightarrow H$ that satisfies
$$\langle T^{\ast}x,y\rangle_{H}=\langle x,Ty\rangle_{H}\quad\forall x,y\in H.$$

In this paper we show that the Cesaro operator
$\mathcal{C}:\ell^{2}\rightarrow\ell^{2}$ is an essentially normal
operator.

\section{preliminaries}

A lot is known about Cesaro operators, most of which come from the
famous paper of Brown, Halmos and Shields published in 1965(See
\cite{BHS}). In that paper Brown, Halmos and Shields\footnote{In
reading that paper I was quite shocked by the fact that although
they cook up everything to conclude that $\mathcal{C}$ is
essentially normal, they do not mention essential normality at
all.} prove that $\mathcal{C}$ is hyponormal, i.e.
$\langle(\mathcal{C}^{\ast}\mathcal{C}-\mathcal{C}\mathcal{C}^{\ast})x,x\rangle_{\ell^{2}}\geq
0$ for all $x\in\ell^{2}$ where $\langle .,.\rangle_{\ell^{2}}$ is
the inner product of $\ell^{2}$ defined as
$$\langle(x_{n})_{n=0}^{+\infty},(y_{n})_{n=0}^{+\infty}\rangle_{\ell^{2}}:=\sum_{n=0}^{+\infty}x_{n}\bar{y_{n}}.$$
In that paper they also compute the spectrum and the point
spectrum of $\mathcal{C}$ and $\mathcal{C}^{\ast}$ as well. They
show that the point spectrum $\sigma_{p}(\mathcal{C})=\emptyset$
of $\mathcal{C}$ is empty whereas the point spectrum
$\sigma_{p}(\mathcal{C}^{\ast})$ of $\mathcal{C}^{\ast}$ is the
open disc
$$\sigma_{p}(\mathcal{C}^{\ast})=\{\lambda\in\mathbb{C}:\mid\lambda-1\mid<1\}$$
with each $\lambda\in\sigma_{p}(\mathcal{C}^{\ast})$ being an
eigenvalue of $\mathcal{C}^{\ast}$ of multiplicity one, i.e. the
dimension of $\ker(\lambda I-\mathcal{C}^{\ast})$ is one for all
$\lambda\in\sigma_{p}(\mathcal{C}^{\ast})$. Also the spectrum
$\sigma(\mathcal{C})$ and $\sigma(\mathcal{C}^{\ast})$ of
$\mathcal{C}$ and $\mathcal{C}^{\ast}$ are equal and that set is
the closed disc
$$\sigma(\mathcal{C})=\sigma(\mathcal{C}^{\ast})=\{\lambda\in\mathbb{C}:\mid\lambda-1\mid\leq 1\}.$$
They also calculate the matrix representations of
$\mathcal{C}^{\ast}\mathcal{C}$ and
$\mathcal{C}\mathcal{C}^{\ast}$ as
$$[\mathcal{C}^{\ast}\mathcal{C}]:=
\begin{pmatrix}
s_{1}&s_{2}&s_{3}&s_{4}&\dots\\
s_{2}&s_{2}&s_{3}&s_{4}&\dots\\
s_{3}&s_{3}&s_{3}&s_{4}&\dots\\
\vdots & \vdots & \vdots &\vdots\\
s_{n}&s_{n}&s_{n}&s_{n}&\dots&s_{n}&s_{n+1}&\dots\\
\vdots & \vdots & \vdots &\vdots
\end{pmatrix}$$
and
$$[\mathcal{C}\mathcal{C}^{\ast}]:=
\begin{pmatrix}
1&\frac{1}{2}&\frac{1}{3}&\frac{1}{4}&\dots\\
\frac{1}{2}&\frac{1}{2}&\frac{1}{3}&\frac{1}{4}&\dots\\
\frac{1}{3}&\frac{1}{3}&\frac{1}{3}&\frac{1}{4}&\dots\\
\vdots & \vdots & \vdots &\vdots\\
\frac{1}{n}&\frac{1}{n}&\frac{1}{n}&\frac{1}{n}&\dots&\frac{1}{n}&\frac{1}{n+1}&\dots\\
\vdots & \vdots & \vdots &\vdots
\end{pmatrix}$$
where
$$s_{n}:=\sum_{k=n}^{+\infty}\frac{1}{k^{2}}.$$
In dealing with the problem of boundedness of $\mathcal{C}$,
Brown, Halmos and Shields use a form of Schur's boundedness test
which I quote below exactly as they formulated it:
\begin{schur}
Let $(X,\Sigma,\mu)$ be a measure space and $k:X\times
X\rightarrow\mathbb{R}$ be a positive measurable function, i.e.
$k(x,y)\geq 0$ for almost every $(x,y)\in X\times X$. Also let
$p:X\rightarrow\mathbb{R}$ be a positive measurable function, i.e.
$p(x)\geq 0$ for almost every $x\in X$. If there are constants
$\alpha$, $\beta\in\mathbb{R}^{+}$ i.e. $\alpha>0$, $\beta>0$ so
that
$$\int_{X}k(x,y)p(y)d\mu(y)\leq\alpha p(x)$$
and
$$\int_{X}k(x,y)p(x)d\mu(x)\leq\beta p(y)$$
are satisfied then the operator $A:L^{2}(X,d\mu)\rightarrow
L^{2}(X,d\mu)$ defined as below
$$(Af)(x):=\int_{X}k(x,y)f(y)d\mu(y)$$
is bounded and the operator norm $\parallel A\parallel$satisfies
$$\parallel A\parallel^{2}\leq\alpha\beta.$$
\end{schur}
  We interpret the Schur's boundedness test above to our situation
  as follows: Here we have $X=\mathbb{N}$ the set of natural
  numbers and $\Sigma=\mathcal{P}(\mathbb{N})$ is the power set of
  $\mathbb{N}$ and $\mu$ is the counting measure. Then if we
  choose $p:\mathbb{N}\rightarrow\mathbb{R}$ to be $p(n)=1$ for
  all $n\in\mathbb{N}$, Schur boundedness test in our setup
  becomes; whenever an operator $A$ on $\ell^{2}$ is given by a
  matrix $[A]:=[a_{i,j}]_{i,j=1}^{+\infty}$ so that
  $$\sup_{n\in\mathbb{N}}\{\sum_{j=1}^{+\infty}\mid a_{n,j}\mid\}:=\alpha<+\infty$$
  and
  $$\sup_{n\in\mathbb{N}}\{\sum_{i=1}^{+\infty}\mid a_{i,n}\mid\}:=\beta<+\infty$$
  are satisfied then $A$ acts boundedly on $\ell^{2}$, i.e.
  $A:\ell^{2}\rightarrow\ell^{2}$ is bounded and the operator norm
  $\parallel A\parallel$ satisfies
  $$\parallel A\parallel^{2}\leq\alpha\beta.$$
  Indeed Schur's boundedness test also gives a compactness
  criterion for operators acting on $\ell^{2}$, i.e. one can show
  that for $[A]:=[a_{i,j}]_{i,j=1}^{+\infty}$ given by the matrix
  $[a_{i,j}]$, if
  $$\lim_{N\rightarrow +\infty}\sup_{n>N}\{\sum_{j=1}^{+\infty}\mid a_{n,j}\mid\}=\lim_{N\rightarrow +\infty}\sup_{n>N}\{\sum_{i=1}^{+\infty}\mid a_{i,n}\mid\}=0$$
  is satisfied then $A:\ell^{2}\rightarrow\ell^{2}$ is a compact operator. To
  show this we look at the truncations of $A$: Let $A_{N}$ be
  given by the matrix $[A_{N}]:=[a_{i.j}^{(N)}]$ so that
  $a_{i,j}^{(N)}:=a_{i,j}$ whenever $i\leq N$ and $j\leq N$. For
  $i>N$ or $j>N$ let $a_{i,j}^{(N)}:=0$. The operator $A_{N}$
  given by such a matrix will have a finite rank and hence is
  compact. So the above condition, by Schur's boundedness test,
  $$\lim_{N\rightarrow +\infty}\sup_{n>N}\{\sum_{j=1}^{+\infty}\mid a_{n,j}\mid\}=\lim_{N\rightarrow +\infty}\sup_{n>N}\{\sum_{i=1}^{+\infty}\mid a_{i,n}\mid\}=0$$
  implies that
  $$\lim_{N\rightarrow +\infty}\parallel A-A_{N}\parallel=0.$$
  Since the ideal of compact operators is norm closed in
  $B(\ell^{2}):=\{A:\ell^{2}\rightarrow\ell^{2}:A\quad\textrm{is linear and
  bounded}\}$, this implies that $A$ is compact as well. We will
  use this argument in proving the compactness of
  $\mathcal{C}^{\ast}\mathcal{C}-\mathcal{C}\mathcal{C}^{\ast}$.
  Let us formulate this as the following corollary of the Schur's
  boundedness test:
  \begin{corollary}
  Let $A$ be a linear operator on $\ell^{2}$ given by the matrix
  $[A]:=[a_{i,j}]$. If the condition
  $$\lim_{N\rightarrow +\infty}\sup_{n>N}\{\sum_{j=1}^{+\infty}\mid a_{n,j}\mid\}=\lim_{N\rightarrow +\infty}\sup_{n>N}\{\sum_{i=1}^{+\infty}\mid a_{i,n}\mid\}=0$$
  is satisfied then $A:\ell^{2}\rightarrow\ell^{2}$ is a compact
  operator.
  \end{corollary}

  We end this preliminary section by pointing out another important
  property of the Cesaro operator, namely its subnormality: An
  operator $T:H\rightarrow H$ on a Hilbert space $H$ over the
  field  of complex numbers $\mathbb{C}$ is called subnormal if
  there is a bigger Hilbert space $\mathcal{H}\supset H$ and a
  normal operator $N:\mathcal{H}\rightarrow\mathcal{H}$,i.e.
  $N^{\ast}N=NN^{\ast}$ so that $H$ is invariant under $N$, i.e.
  $N(H)\subseteq H$ and the restriction of $N$ to $H$ coincides with $T$,
  i.e. $N|_{H}=T$. Kriete and Trutt proved that the Cesaro
  operator $\mathcal{C}$ is subnormal(see \cite{KT}) in 1971.

  \section{The Main Result}

  We now prove the main result of the paper:
\begin{thma}
The Cesaro operator $\mathcal{C}:\ell^{2}\rightarrow\ell^{2}$ is
essentially normal, i.e. the self commutator
$[\mathcal{C}^{\ast},\mathcal{C}]:=\mathcal{C}^{\ast}\mathcal{C}-\mathcal{C}\mathcal{C}^{\ast}$
is a compact operator on $\ell^{2}$.
\end{thma}
\begin{proof}
We prove the compactness of
$\mathcal{C}^{\ast}\mathcal{C}-\mathcal{C}\mathcal{C}^{\ast}$ by
looking at its matrix representation which was calculated by
Brown, Halmos and Shields as follows:
$$[\mathcal{C}^{\ast}\mathcal{C}-\mathcal{C}\mathcal{C}^{\ast}]:=
\begin{pmatrix}
s_{1}-1&s_{2}-\frac{1}{2}&s_{3}-\frac{1}{3}&s_{4}-\frac{1}{4}&\dots\\
s_{2}-\frac{1}{2}&s_{2}-\frac{1}{2}&s_{3}-\frac{1}{3}&s_{4}-\frac{1}{4}&\dots\\
s_{3}-\frac{1}{3}&s_{3}-\frac{1}{3}&s_{3}-\frac{1}{3}&s_{4}-\frac{1}{4}&\dots\\
\vdots & \vdots & \vdots &\vdots\\
s_{n}-\frac{1}{n}&s_{n}-\frac{1}{n}&s_{n}-\frac{1}{n}&s_{n}-\frac{1}{n}&\dots&s_{n}-\frac{1}{n}&s_{n+1}-\frac{1}{n+1}&\dots\\
\vdots & \vdots & \vdots &\vdots
\end{pmatrix}
$$
 where
$$s_{n}:=\sum_{k=n}^{+\infty}\frac{1}{k^{2}}. $$

 The self commutator is self adjoint so we have for
$[\mathcal{C}^{\ast}\mathcal{C}-\mathcal{C}\mathcal{C}^{\ast}]:=[a_{i,j}]$,
$$\sum_{j=1}^{+\infty}\mid a_{n,j}\mid=\sum_{i=1}^{+\infty}\mid a_{i,n}\mid\quad\forall n\in\mathbb{N}\setminus\{0\}$$
Hence it is enough to show that
$$\lim_{N\rightarrow +\infty}\sup_{n>N}\{\sum_{j=1}^{+\infty}\mid a_{n,j}\mid\}=0.$$
We observe that for $j\leq n$ we have $a_{n,j}=s_{n}-\frac{1}{n}$
and for $j>n$ we have $a_{n,j}=s_{j}-\frac{1}{j}$. So this gives
$$\sum_{j=1}^{+\infty}\mid a_{n,j}\mid=n(s_{n}-\frac{1}{n})+\sum_{j=n+1}^{+\infty}s_{j}-\frac{1}{j}\quad\forall n\in\mathbb{N}\setminus\{0\}$$
Having
$$\lim_{n\rightarrow +\infty}\sum_{j=1}^{+\infty}\mid a_{n,j}\mid=0$$
will imply
$$\lim_{N\rightarrow +\infty}\sup_{n>N}\{\sum_{j=1}^{+\infty}\mid a_{n,j}\mid\}=0$$
hence it is enough to prove that
$$\lim_{n\rightarrow +\infty}\sum_{j=1}^{+\infty}\mid a_{n,j}\mid=0.$$
Let us first look at the expression $s_{n}-\frac{1}{n}$: consider
$$\frac{1}{n}=\int_{n}^{+\infty}\frac{dx}{x^{2}}=\sum_{k=n}^{+\infty}\int_{k}^{k+1}\frac{dx}{x^{2}}.$$
Using this we have
$$s_{n}-\frac{1}{n}=\sum_{k=n}^{+\infty}\frac{1}{k^{2}}-\sum_{k=n}^{+\infty}\int_{k}^{k+1}\frac{dx}{x^{2}}=\sum_{k=n}^{+\infty}\bigg(\frac{1}{k^{2}}-\int_{k}^{k+1}\frac{dx}{x^{2}}\bigg)$$
$$=\sum_{k=n}^{+\infty}\frac{1}{k^{2}}-(\frac{-1}{x})|_{x=k}^{x=k+1}=\sum_{k=n}^{+\infty}\frac{1}{k^{2}}-(\frac{1}{k}-\frac{1}{k+1})=\sum_{k=n}^{+\infty}\frac{1}{k^{2}}-\frac{1}{k(k+1)}$$
$$=\sum_{k=n}^{+\infty}\frac{1}{k^{2}(k+1)}$$
Since $\frac{1}{k+1}<\frac{1}{k}$, we have
$$s_{n}-\frac{1}{n}<\sum_{k=n}^{+\infty}\frac{1}{k^{3}}\quad\forall n\in\mathbb{N}\setminus\{0\}.$$
Since $f(x)=\frac{1}{x^{3}}$ is a decreasing function on positive
real numbers we have
$$\sum_{k=n}^{+\infty}\frac{1}{k^{3}}<\int_{n-1}^{+\infty}\frac{dx}{x^{3}}=\frac{1}{2(n-1)^{2}}\quad\forall n\in\mathbb{N},n>1$$
So we have
$$s_{n}-\frac{1}{n}<\frac{1}{2(n-1)^{2}}\quad\forall n\in\mathbb{N},n>1$$
and this implies that
$$\sum_{j=1}^{+\infty}\mid a_{n,j}\mid=n(s_{n}-\frac{1}{n})+\sum_{j=n+1}^{+\infty}s_{j}-\frac{1}{j}<\frac{n}{2(n-1)^{2}}+\frac{1}{2}\sum_{j=n}^{+\infty}\frac{1}{j^{2}}<\frac{n}{2(n-1)^{2}}+\frac{1}{2(n-1)}$$
$$=\frac{2n-1}{2(n-1)^{2}}$$
So this estimate implies that
$$\lim_{n\rightarrow +\infty}\sum_{j=1}^{+\infty}\mid a_{n,j}\mid=0$$
and hence, by the Corollary 1 of the preliminaries section, the
operator
$\mathcal{C}^{\ast}\mathcal{C}-\mathcal{C}\mathcal{C}^{\ast}$ is
compact.
\end{proof}

\begin{center}
\textbf{Acknowledgements}
\end{center}
The author is indebted to Prof. Aristomenis Siskakis of Aristotle
University of Thessaloniki for pointing out the main result of the
paper.

\end{document}